\documentclass[12pt]{amsart}
\usepackage[margin=0.6in]{geometry}

\usepackage[english]{babel}
\usepackage[alphabetic]{amsrefs}
\usepackage{amsmath,amssymb,amsfonts,amsthm,enumerate}
\usepackage{url}
\usepackage{graphicx}
\usepackage{wrapfig}

\numberwithin{equation}{section}

\newtheorem{theorem}{Theorem}[section]

\newcommand{\R}{\mathbb{R}}

\renewcommand{\varpi}{\omega}

\renewcommand{\le}{\leqslant}

\renewcommand{\ge}{\geqslant}

\renewcommand{\epsilon}{\varepsilon}

\title[Free boundary problems]{(Non)local and (non)linear
free boundary problems}

\author{Serena Dipierro and Enrico Valdinoci}

\address{{\em Serena Dipierro:} School of Mathematics and Statistics,
University of Melbourne, 813 Swanston St, 
Parkville VIC 3010, Australia}
\email{sdipierro@unimelb.edu.au}

\address{{\em Enrico Valdinoci:} School of Mathematics and Statistics,
University of Melbourne, 813 
Swanston St, Parkville VIC 3010, Australia,
Dipartimento di Matematica, Universit\`a degli studi di Milano,
Via Saldini 50, 20133 Milan, Italy}


\begin{document}

\begin{abstract}
We discuss some recent developments in the theory
of free boundary problems, as obtained in a series of
papers in collaboration with L.~Caffarelli,
A.~Karakhanyan and O.~Savin. 

The main feature of these new free boundary problems
is that they deeply take into account nonlinear energy superpositions
and possibly nonlocal functionals.

The nonlocal parameter interpolates between
volume and perimeter functionals, and so it can be seen
as a fractional counterpart of classical
free boundary problems, in which the bulk energy
presents nonlocal aspects.

The nonlinear term in the energy superposition
takes into account the possibility of modeling different
regimes in terms of different energy levels and provides
a lack of scale invariance, which in turn may cause
a structural instability of minimizers that may vary from
one scale to another.
\end{abstract}
\maketitle
\tableofcontents

\hfill
\begin{quote}
\texttt{{
``The shape of a snowflake as developed by 
the freezing process is governed by the distribution
of temperature within the water particles 
and by the law of conservation of energy along the boundary. 
This shape is not known in
ad\-vance, and its surface 
is a classical example of a free boundary formed by nature.''\\
(Avner Friedman, \cite{MR1776102})}}
\end{quote}

\bigskip

\section{Introduction: classical
linear free boundary problems}

\subsection{Linear local plus volume free boundary problems}

A free boundary problem is 
a mathematical problem with (at least) two unknowns,
typically a function and a domain.
The domain may depend on the function itself,
which is usually determined by
a partial differential equation, whose boundary conditions
depend in turn on the domain.
As a typical example, one may consider the problem
of finding a
function~$u$, whose values are prescribed
equal to some function~$u_o$
at the boundary of a domain~$\Omega\subset\R^n$,
which is harmonic outside its zero level set
and with constant normal derivative along its zero level set,
namely
\begin{equation}\label{FB} \left\{\begin{matrix}
& u=u_o& {\mbox{ on }}\partial\Omega,\\
& \Delta u=0 &{\mbox{ in }} \Omega\setminus\{u=0\},\\
& |\nabla u^+|^2 -|\nabla u^-|^2=\Lambda>0 &
{\mbox{ on }}\Omega\cap \big( \partial\{u>0\}\cup
\partial\{u<0\}
\big).
\end{matrix}
\right.\end{equation}
Here, we adopted the standard notation~$u^+:=\max\{u,0\}$
and~$u^-:= \max\{ -u,0\}$.
The set~$\Omega\cap \big( \partial\{u>0\}\cup
\partial\{u<0\}
\big)$ is called the free boundary.
In problem~\eqref{FB}, the interesting unknowns are both~$u$
and its positivity set~$\Omega\cap\{u>0\}$:
of course, if we know~$u$ to start with, we also know
its positivity set and, viceversa,
if we know its positivity set we can 
solve the associated partial differential equation
to find the harmonic function with
prescribed boundary data. Nevertheless,
it is convenient to look at problem~\eqref{FB}
as a whole,
in which a PDE with boundary conditions
is supplemented by
an additional condition at the free boundary.
Such condition is often motivated by physical applications:
for instance,
the normal derivative condition in~\eqref{FB}
is a consequence of the Bernoulli's Law in the framework
of ideal fluid jets, see e.g. the 
discussion in~\cite{2017arXiv170107897D}.
See also~\cite{MR1776102}
for an exhaustive introduction to free boundary problems
also in view of concrete applications.
\medskip

Following the classical works of~\cite{MR618549}
and~\cite{MR732100}, it is convenient to
set problem~\eqref{FB} into a somehow 
more convenient variational framework.
Namely, one considers the energy functional
\begin{equation}\label{ACF}
\int_\Omega |\nabla u(x)|^2\,dx+
\lambda_+ {\mathcal{L}}^n\big(\Omega_+(u)\big)
+\lambda_-{\mathcal{L}}^n\big(\Omega_-(u)\big),
\end{equation}
with~$\Lambda:=\lambda_+-\lambda_->0$,
$\Omega_+(u) := \Omega\cap \{u>0\}$
and
$\Omega_-(u) :=\Omega\cap \{u<0\}$,
being~${\mathcal{L}}^n$ the Lebesgue measure.
One looks at minimizers of~\eqref{ACF}
with~$u=u_o$ along~$\partial\Omega$ (say, in the trace sense).
It is proved in Theorem~2.2
of~\cite{MR732100} that
minimizers are harmonic outside their zero set.
Also, in Theorem~2.4 of~\cite{MR732100} it is shown that
\begin{equation}\label{BERN-AC}
|\nabla u^+|^2 -|\nabla u^-|^2 =\Lambda\end{equation}
along the free boundary (in an appropriate weak
sense and provided that the zero
set has zero Lebesgue measure). Therefore,
the minimization of~\eqref{ACF} provides
a concise and
useful variational tool to set problem~\eqref{FB}
into an appropriate functional setting.
In addition, the energy framework in~\eqref{ACF}
makes a number of techniques from harmonic analysis,
elliptic PDEs and geometric measure theory
available to study free boundary problems.
In view of this, it was possible in~\cite{MR732100}
to establish a series of rigidity and regularity results
for minimizers, including
{\tt ``nondegeneracy 
theorems, the Lipschitz continuity of the solution, 
and some properties of the free boundary; for~$n = 2$ 
the free boundary is proved to be continuously differentiable.

A new and rather powerful tool introduced in this paper 
is the monotonicity formula''}, which says that if~$u$
is a minimizer and~$0$ belongs to the free boundary, then the
quantity
\begin{equation}\label{0u854udyy6gyitirghioer1203498whc} \frac{1}{r^4} \int_{B_r} \frac{|\nabla u^+(x)|^2}{|x|^{n-2}}\,dx
\,\int_{B_r} \frac{|\nabla u^-(x)|^2}{|x|^{n-2}}\,dx\end{equation}
is monotone increasing in~$r$.\medskip

\begin{wrapfigure}{L}{9cm}
    \centering
    \includegraphics[width=9cm]{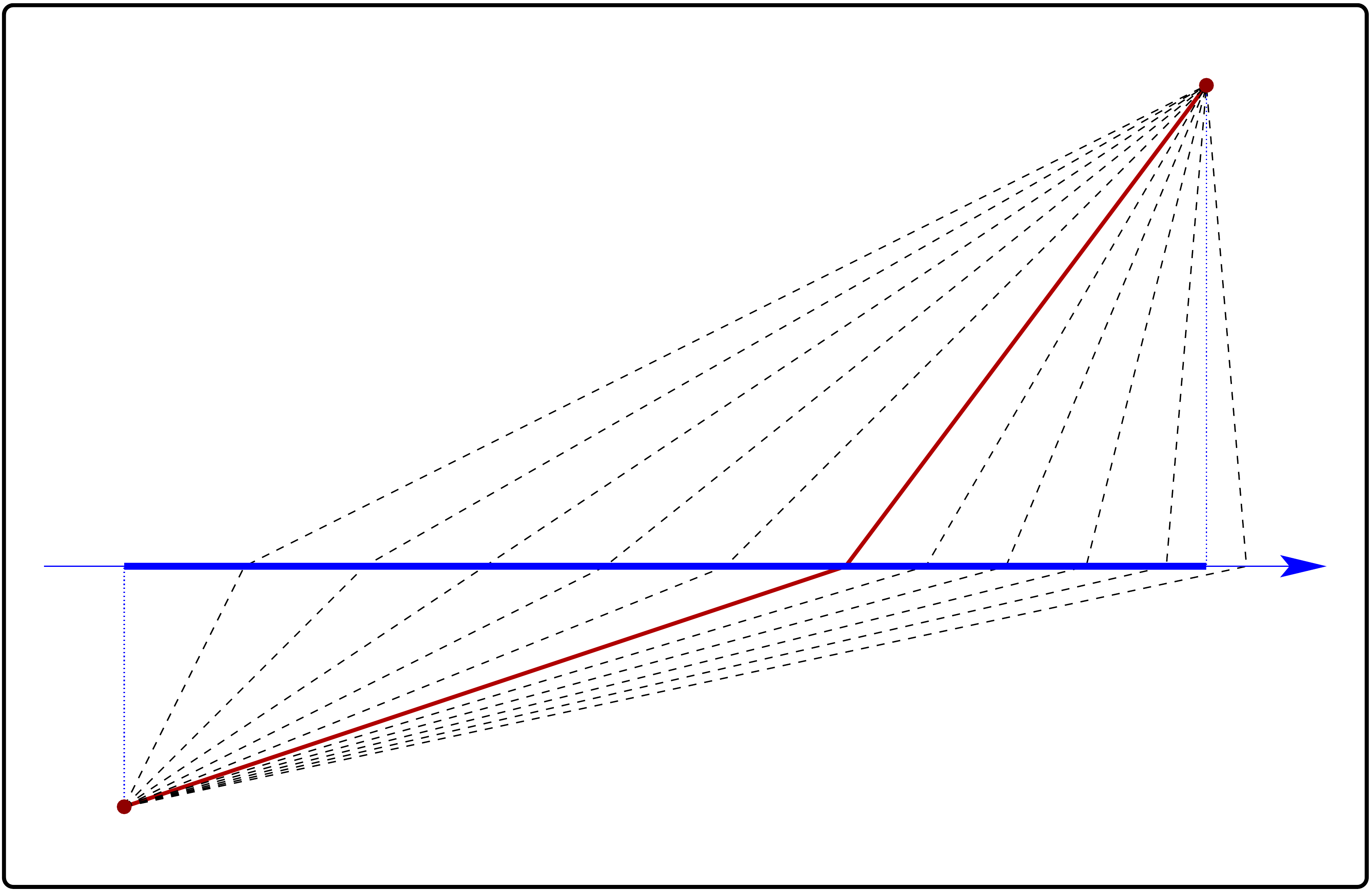}
\end{wrapfigure}

A heuristic simplification of this type of problems can be easily
discussed by looking at the one-dimensional case. In this
setting, the harmonicity of the minimizers outside the zero set
gives that the graph of the
solution,
with boundary data prescribed outside an interval,
is obtained by glueing two segments, which join the boundary data
to an intersection with the horizontal axis. The position of this
intersection point is then chosen to minimize the total energy functional,
and this imposes
a prescription between the inward and outward slope of the minimizing function
at its zero set, see the figure to the side.
In this sense, the Bernoulli condition at the free boundary translates into
a prescription of the slope of the minimizer, hence on a Lipschitz and nondegeneracy
property (of course, the higher dimensional case becomes highly nontrivial,
but the spirit of the proofs is often to ``try to reduce the picture to the one-dimensional case'',
by blow-up and classification methods).\medskip

The case in which the boundary datum~$u_o$ is nonnegative
(or equivalently, by the maximum principle, the case
in which the minimizer is nonnegative) is called the ``one-phase case'',
the phase being that represented by the positive set of
the minimizer;
in this jargon, the ``two-phase case'' corresponds to the case
in which both the positive and the negative set
of the minimizer are nontrivial (for some reason, no phase
is associated to the zero level set in this notation).

It is worth remarking that in the one-phase case~$\Omega^-(u)$
is trivial and so the functional in~\eqref{ACF}
reduces to that studied in~\cite{MR618549}, namely
\begin{equation}\label{AC}
\int_\Omega |\nabla u(x)|^2\,dx+
\lambda_+ {\mathcal{L}}^n\big(\Omega_+(u)\big)
.\end{equation}
It is also interesting to observe that the energy functional
in~\eqref{AC} (and correspondingly that in~\eqref{ACF},
with minor modifications in the discussion) can be seen as
the ``superposition'' of two energy terms:
namely a Dirichlet energy term, given by the seminorm in
the Sobolev space~$H^1(\Omega)$, and a bulk energy
term of volume type, given by the measure of the positivity set.
{F}rom the perspective of this paper, the superposition
of Dirichlet and bulk energy in~\eqref{AC}
is ``linear'', i.e. the total energy is linear in each of
the single energy contributions.
\medskip

\subsection{Linear local plus perimeter free boundary problems}

A more recent variant of~\eqref{AC} has been introduced in~\cite{MR1808651}
and it deals with the case in which the volume term is replaced
by a perimeter term, namely one takes into account energy
functionals of the form
\begin{equation}\label{SALSA}
\int_\Omega |\nabla u(x)|^2\,dx+
{\rm Per}\big(\{u>0\},\Omega\big),
\end{equation}
where~${\rm Per}$ denotes the standard perimeter functional
(see e.g.~\cite{MR0638362}).
As a matter of fact, the right functional setting
for~\eqref{SALSA} is slightly more involved, since
it must take into account separately the function~$u$
and a set~$E$ on which~$u\ge0$ (rather than just
the na\"{\i}ve
positivity set of~$u$), but, for the sake of simplicity
we will not indulge in
this technical detail in this paper.

It is shown in \cite{MR1808651} that
the minimizers of~\eqref{SALSA} enjoy suitable regularity
properties: in particular, the authors show that the Bernoulli condition
along the free boundary reads
\begin{equation}\label{BERN-SALSA}
|\nabla u^+|^2 -|\nabla u^-|^2 =H_u,\end{equation}
where~$H_u$ is the mean curvature of the zero level set of~$u$
(one can compare
this with~\eqref{BERN-AC});
also, the ``natural scaling'' of the problem is related to the 
H\"older power~$\frac12$, and it is proved 
in~\cite{MR1808651} that minimizers are H\"older -- and, in fact,
Lipschitz, thus ``beating the natural scale of the problem''.

With this, one obtains that the gradient terms in~\eqref{BERN-SALSA}
act as a ``lower order term'',
thus the free boundary regularity is reduced to that of minimal
surfaces.\medskip

It is also worth repeating that the functional in~\eqref{SALSA}
is also the linear superposition of two energies,
namely the Dirichlet energy given by the seminorm
in~$H^1(\Omega)$ and the perimeter functional.\medskip

The goal of this paper is to present a series of
very recent results in which
one takes into account ``nonlinear superpositions''
of energies. In addition, we will consider the cases
in which one of these energies, or both, are ``nonlocal'',
in a sense that will be clear by the context in the sequel.

\section{Nonlinear local plus volume free boundary problems}

In \cite{2016arXiv161100412D} a new type of free boundary
problems has been introduced, in which
the two energy terms in the total functional are superposed
in a nonlinear way. That is, one of the two energy terms
``counts'' more (or less) than the other one in suitable
energy ranges. These problems aim to model cases
in which the Dirichlet or the bulk energies become predominant
when different scales are considered or when
different energy levels are taken into account.\medskip

In particular, in~\cite{2016arXiv161100412D} one considers the
energy functional
\begin{equation}\label{6283ejdefw}
\int_\Omega |\nabla u(x)|^2\,dx+
\Phi\big(\lambda_-\Omega_-(u),\,\lambda_+\Omega_+(u)\big).
\end{equation}
Here, $\Phi$ is a ``nonlinearity'',
which is supposed to be sufficiently smooth, normalized
in such a way~$\Phi(0,0)=0$, and satisfying the monotonicity
property
$$ \Phi_0'(r)\ge0,$$
where
$$ \Phi_0(r):=\Phi\left( \lambda_-\left( {\mathcal{L}}^n(\Omega)
-\frac{r}{\lambda_+}\right),\,r\right) .$$
This notation is also useful to write
the volume energy in~\eqref{6283ejdefw}
in terms only of the positive bulk, namely to reduce~\eqref{6283ejdefw}
to the energy functional
\begin{equation*}
\int_\Omega |\nabla u(x)|^2\,dx+
\Phi_0\big(\lambda_+\Omega_+(u)\big).
\end{equation*}
We remark that when~$\Phi(r_1,r_2):=r_1+r_2$
the functional in~\eqref{6283ejdefw}
reduces to the classical one in~\eqref{ACF}.

In addition, when~$\Phi(r_1,r_2)\sim r_2^{\frac{n-1}n}$
or~$\Phi(r_1,r_2)\sim (r_1+r_2)^{\frac{n-1}n}$
the functional in~\eqref{6283ejdefw}
provides a scaling between volume and perimeter
which is naturally
related to the isoperimetric inequality.\medskip

We remark that when~$\Phi_0$ is concave, then
the minimizers
of the nonlinear energy functional in~\eqref{6283ejdefw}
boil down to those in \eqref{AC}
(see the computation after Theorem~1.5
in~\cite{2016arXiv161100412D}).

Nevertheless, the general case provides additional source of
difficulties, due to the instability of the minimizers
and the lack of scale invariance of the problem. More explicitly,
the nonlinear structure of the problem produces the feature that
a minimizer in a given domain is not necessarily a minimizer in a
smaller subdomain, see e.g. the example in Theorem~1.1
of~\cite{2015arXiv151203043D}.\medskip

The Bernoulli type free boundary condition
is also more complicated in the nonlinear case, and it takes
into account global properties of the minimizer itself,
while, for instance, the condition in~\eqref{BERN-AC} is independent of the minimizer
and only depends on structural constants, and the one in~\eqref{BERN-SALSA}
depends only on a local feature of the minimizer, such as the mean curvature of its free boundary.
Without going into the details of the weak formulation of the meaning of
the Bernoulli condition at irregular free boundary points, we mention here that,
as proved in~\cite{2016arXiv161100412D}, at regular free boundary points,
the minimizers satisfy
\begin{equation}\label{self cond}
|\nabla u^+|^2 -|\nabla u^-|^2=
\lambda_+ \partial_{r_2} \Phi\big(\lambda_-\Omega_-(u),\,\lambda_+\Omega_+(u)\big)
-\lambda_-\partial_{r_1}
\Phi\big(\lambda_-\Omega_-(u),\,\lambda_+\Omega_+(u)\big).
\end{equation}
Notice that~\eqref{self cond}, in general, takes into account a global property of a minimizer,
namely its positive and negative masses~$\Omega_-(u)$
and~$\Omega_+(u)$, nevertheless if~$\Phi(r_1,r_2)=r_1+r_2$ then~\eqref{self cond}
reduces to the classical one in~\eqref{BERN-AC}.\medskip

The nonlinear structure of
the energy functional in~\eqref{6283ejdefw} is also a source of serious issues
towards regularity, even in low dimension. For example, in Section~4 of~\cite{2016arXiv161100412D}
it is proved that there exists a nonlinearity~$\Phi$ such that
the function~$u(x)=u(x_1,x_2):=x_1x_2$
is a minimizer in the unit ball of the plane.\medskip

On the other hand, under a monotonicity assumption on the nonlinearity
(as stated in~\eqref{MONP} here below),
a series of positive regularity results has been established in 
Theorems~1.1, 1.3, 1.4, 1.6 and~1.7
in~\cite{2016arXiv161100412D}, which can be summarized as follows:

\begin{theorem}\label{09asudy98r389rfgas}
Let $u$ be a minimizer in~$\Omega\subset\R^n$
for the functional in~\eqref{6283ejdefw} with
\begin{equation}\label{MONP} \inf \Phi_0'>0.\end{equation}
Suppose that~$0\in\partial\Omega_+(u)$
and~${\mathcal{L}}^n(\Omega_+(u))>0$.

Then:
\begin{itemize}
\item (Regularity of the minimizers). $u$ is locally Lipschitz continuous in~$\Omega$.
\item (Nondegeneracy of the minimizers). For small~$r>0$, there exists~$c>0$
such that 
$$ \int_{B_r\cap \{u>0\}} u^2(x)\,dx\ge c r^{n+2}.$$
\item (Regularity of the positive phase). $\Delta u^+$ is a Radon measure and,
for small~$r>0$, there exists~$c>0$
such that 
$$ cr^{n-1}\le\int_{B_r} \Delta u^+(x)\,
dx\le\frac{1}{r} \int_{B_{2r}} |\nabla u^+(x)|\,dx.$$
\item (Regularity of the free boundary). $\partial\Omega^+(u)$ has locally
finite $(n-1)$-dimensional Haussdorff measure. Also
the points of the reduced boundary have full $(n-1)$-dimensional Haussdorff measure
in~$\partial\Omega^+(u)$.
\item (Regularity of the free boundary in the plane). If~$n=2$, then~$\partial\Omega^+(u)$
is a continuously differentiable curve.
\end{itemize}\end{theorem}

It is worth to point out that minimizers of nonlinear free boundary problems
reduce to minimizers of classical free boundary problems in the blow-up limit.
Namely, under the assumptions of Theorem~\ref{09asudy98r389rfgas},
one can consider an infinitesimal sequence~$\rho_k$ and define the blow-up sequence
\begin{equation}\label{BUP} u_k(x):=\frac{u(\rho_k x)}{\rho_k}\end{equation}
and obtain that~$u_k$ converges to some~$u_0$ locally uniformly in~$\R^n$ and that~$u_0$
is a minimizer of the functional in~\eqref{AC}. See~\cite{2016arXiv161100412D}
for further details.

\section{Linear nonlocal plus volume free boundary problems}

In~\cite{MR2677613}, \cite{MR2926238} and~\cite{MR3234971}
a nonlocal variation of the classical energy functional in~\eqref{AC}
has been considered. Such functional can be written as
\begin{equation}\label{FBS}
\int_{\Omega^+} t^{1-2s} |\nabla U(x,t)|^2\,dx\,dt+
{\mathcal{L}}^n \big( \{u>0\}\cap\Omega_0\big),
\end{equation}
where~$s\in(0,1)$,
$(x,t)\in\R^{n+1}_+:=\R^n\times(0,+\infty)$,
$u(x)=U(x,0)$, $\Omega\subset\R^{n+1}$, $\Omega^+:=\Omega\cap\{t>0\}$
and~$\Omega_0:=\Omega\cap\{t=0\}$.
Several regularity results
are shown in~\cite{MR2677613}, \cite{MR2926238} and~\cite{MR3234971},
such as
optimal regularity, nondegeneracy
and smoothness of the free boundary.

In terms of equations and free boundary conditions,
the minimization of~\eqref{FBS} is related (up to dimensional
constants that we omit for the sake of simplicity) to the fractional Laplacian,
and can be formally written as
\begin{equation} \label{BCS}\left\{\begin{matrix}
& (-\Delta)^s u=0 &{\mbox{ in }} \Omega_0\cap\{u>0\},\\
& \displaystyle\lim_{y\to x}
\frac{u(y)}{\big((y-x)\cdot\nu(x)\big)^s} =1&
{\mbox{ if }}x\in\Omega_0\cap \big( \partial\{u>0\}
\big),
\end{matrix}
\right.\end{equation}
see formula~(1.2) in~\cite{MR2677613} and the comments
therein for further details about
this. In~\eqref{BCS}, $\nu(x)$ denotes the unit normal at a point~$x$
of the free boundary, pointing towards the positivity set
of~$u$. Also, we used the standard notation for the fractional
Laplacian, namely
\begin{equation}\label{DEF:0123rugrfhijfoo98utytg}
(-\Delta)^s u(x):=\int_{\R^n} \frac{2u(x)-u(x+y)-u(x-y)}{|y|^{n+2s}}\,dy
.\end{equation}
The reader may compare \eqref{BCS} with \eqref{FB}.
Also, the extended function~$U$ is related to the trace function~$u$
by a singular or degenerate problems with weights, namely
\begin{equation*}
\left\{\begin{matrix}
& {\rm div}\,(t^{1-2s}\nabla U(x,t))=0 &{\mbox{ in }} 
\R^{n+1}_+,\\
& U(\cdot,0)=u(\cdot)&
{\mbox{ in }}\R^n,
\end{matrix}
\right.\end{equation*}
see e.g. formula~(1.5) in~\cite{MR2677613}.\medskip

It is interesting to point out that
the functional in~\eqref{FBS} is the sum (thus the linear superposition)
of a nonlocal Dirichlet energy and a volume term, whence
the title of this section.

\section{Linear local plus nonlocal free boundary problems}

In \cite{MR3390089} the free boundary problem
arising from the superposition of a classical Dirichlet energy
and a fractional perimeter has been introduced.
In this case, one considers the energy functional
\begin{equation}\label{PSIGMA}
\int_\Omega |\nabla u(x)|^2\,dx+
{\rm Per}_\sigma\big( \{u>0\},\Omega\big).
\end{equation}
Here~$\sigma\in(0,1)$ is a fractional parameter and~${\rm Per}_\sigma$
denotes the fractional perimeter introduced in~\cite{MR2675483}.
Namely, one defines the fractional interaction of two disjoint subsets~$E$ and~$F$
of~$\R^n$ as
$$ {\mathcal{I}}_{\sigma}(E,F):=
\sigma\,(1-\sigma)\, 
\iint_{E\times F}\frac{dx\,dy}{|x-y|^{n+\sigma}}.$$
Then,
fixed a reference domain~$\Omega$ and a set~$E\subseteq\R^n$
one defines the $\sigma$-perimeter of~$E$ in~$\Omega$ as the sum of all
the contributions in the fractional interaction between~$E$ and its complement~$E^c$
which involve
the domain~$\Omega$, that is
$$ {\rm Per}_\sigma(E,\Omega):=
{\mathcal{I}}_{\sigma}(E\cap\Omega, E^c\cap\Omega)+
{\mathcal{I}}_{\sigma}(E\cap\Omega, E^c\cap \Omega^c)+
{\mathcal{I}}_{\sigma}(E\cap\Omega^c,E^c\cap\Omega).$$
The study of the fractional minimal surfaces (i.e. of the minimizers
of the fractional perimeter) is a very interesting topic of investigation in itself, which
offers a great number of extremely challenging and important open problems
(we refer to~\cite{GP123456789} for a recent review on this topic).

In particular, we recall (see e.g.~\cites{MR3586796, MR1942130, MR2782803}
and~\cites{MR1940355, MR3007726})
that the fractional perimeter recovers the classical one as~$\sigma\nearrow1$
and, in some sense, it recovers the Lebesgue measure as~$\sigma\searrow0$.
In this sense, the functional in~\eqref{PSIGMA} may be seen as a fractional interpolation
of those in~\eqref{AC} and~\eqref{SALSA}, which are extrapolated from~\eqref{PSIGMA}
in the limit as~$\sigma\searrow0$ and as~$\sigma\nearrow1$, respectively.\medskip

The Bernoulli condition along the free boundary for
the functional in~\eqref{PSIGMA} is now of nonlocal type.
Namely, one can consider the fractional mean curvature
of a set~$E$
at a point~$x\in\partial E$, that is
$$ H_\sigma(x,E):=
\sigma\,(1-\sigma)\,
\int_{\R^n} \frac{\chi_{E^c}(y)-\chi_E(y)}{|x-y|^{n+\sigma}}\,dy.$$
This quantity plays a special role concerning the interior
and boundary behaviors of fractional minimal surfaces
(see e.g.~\cite{MR2675483} and~\cite{MR3596708}).
In addition
(when~$E$ is bounded and sufficiently regular) it
recovers the classical mean curvature as~$\sigma\nearrow1$
and approaches a constant as~$\sigma\searrow0$ (see the
appendices in~\cite{GP123456789}; see also~\cite{MR3230079}
for the basic properties of~$H_\sigma(x,E)$).
\medskip

In our framework, the fractional mean curvature plays a crucial role
for the free boundary condition, which, for
the functional in~\eqref{PSIGMA}, reads
\begin{equation}\label{HGA:0}
|\nabla u^+(x)|^2 -|\nabla u^-(x)|^2 =H_\sigma(x,E),
\end{equation}
being~$u$ a minimizer with~$u\ge0$ in~$E$.
We stress that~\eqref{HGA:0} reproduces
the classical Bernoulli conditions
in~\eqref{BERN-AC} 
and~\eqref{BERN-SALSA}
when~$\sigma\searrow0$ and~$\sigma\nearrow1$, respectively.\medskip

Furthermore, one can obtain the following
regularity results, due to Theorems~1.1 and~1.2  
in~\cite{MR3390089} (once again, we omit the details related
to the notion of minimizing pair for simplicity):

\begin{theorem}\label{SIG:THM}
Let $u$ be a minimizer for~\eqref{PSIGMA}
in~$\Omega\subset\R^n$
with~$u\ge0$ in~$E$ and~$0\in\partial E$.

Then:
\begin{itemize}
\item (Regularity of the minimizers). $u$ is locally H\"older continuous,
with H\"older exponent~$1-\frac\sigma2$.
\item (Density estimates for the free boundary). For
small~$r>0$, 
$$ \min\big\{ {\mathcal{L}}^n (E\cap B_r),\;
{\mathcal{L}}^n (E^c\cap B_r)\big\}\ge cr^n,$$
for some~$c>0$.
\item (Regularity of the free boundary). 
$\partial E$ is a smooth hypersurface 
outside a small singular set of
Hausdorff dimension at most~$n-3$.
\item (Regularity of the free boundary in the plane). If~$n=2$, 
then~$\partial E$
is a smooth hypersurface.
\end{itemize}\end{theorem}

We think that it is a very important open problem
to detect the optimal regularity for the minimizers
and their free boundaries in the setting of Theorem~\ref{SIG:THM}.
Indeed, it would be very desirable to establish full regularity
of the free boundary in higher dimension and/or to provide
counterexample (this problem is likely to be related
to the regularity of the minimizers
or almost minimizers of the fractional perimeter).

Moreover, the optimal regularity for~$u$ in the setting
of Theorem~\ref{SIG:THM} is not known. We point out that the H\"older
exponent~$1-\frac\sigma2$ in Theorem~\ref{SIG:THM} recovers
the Lipschitz exponent in~\cite{MR618549} as~$\sigma\searrow0$,
but not the Lipschitz exponent in~\cite{MR1808651}
as~$\sigma\nearrow1$, therefore it is very plausible that
the regularity theory of Theorem~\ref{SIG:THM} can be improved.
Such improvement, however, cannot be based simply on scaling properties,
since the natural scale of problem~\eqref{PSIGMA}
relies on the exponent~$1-\frac\sigma2$, and the blow-up sequence
in this case is given by
\[ u_k(x):=\frac{u(\rho_k x)}{\rho_k^{1-\frac\sigma2}},\]
to be compared with~\eqref{BUP} -- hence an improvement of regularity
needs in this case
to ``beat'' the natural scaling of the problem.

\section{Linear nonlocal plus nonlocal free boundary problems}

In~\cites{MR3427047, POIPOIPOI}
the following energy functional is taken into consideration:
\begin{equation}\label{9hiocdvjkb2toi3rg2fkfdjhsgddksjaqqd}
{\mathcal{G}}_s (u,\Omega)+
{\rm Per}_\sigma\big( \{u>0\},\Omega\big).
\end{equation}
In this setting, $s$ and~$\sigma$ are fractional
parameter, ranging in the interval~$(0,1)$,
and~${\mathcal{G}}_s$ is a Dirichlet energy form
related to the Gagliardo seminorm, namely
\begin{equation}\label{asdfgv98dguiuiefgu923igfuweioj} {\mathcal{G}}_s (u,\Omega):=\iint_{{\mathcal{Q}}_\Omega}\frac{
|u(x)-u(y)|^2
}{|x-y|^{n+2s}}\,dx\,dy,\end{equation}
with
$$ {\mathcal{Q}}_\Omega:=
(\Omega\times\Omega)\cup (\Omega\times\Omega^c)\cup(\Omega^c\times\Omega).$$
Differently from the classical Dirichlet energy, the
Gagliardo seminorm in~\eqref{asdfgv98dguiuiefgu923igfuweioj}
takes into account a ``weighted'' and ``global'' oscillation property
of the function and its minimization is related to the fractional Laplacian operator.
This modification creates a series of conceptual differences (and also
a great number of technical difficulties): it is indeed worth recalling
that $s$-harmonic functions (i.e. functions with vanishing fractional
Laplacian) behave very differently from the classical harmonic
functions, possess a much richer ``zoology''
and much less local rigidity results. In particular, as shown in~\cite{MR3626547},
{\em any} given function (independently of its geometric properties)
can be locally approximated by a $s$-harmonic function, and this
fact is in sharp contrast with the rigidity features exhibited
by classical harmonic functions.\medskip

The ``abundance'' of $s$-harmonic functions and their lack of regularity
near the boundary provide a series of important obstacles towards general
regularity theories for fractional problems
and the $s$-harmonic replacement problem provides a series of interesting
features, as investigated in~\cites{MR3320130, MR3427047, POIPOIPOI}.\medskip

In~\cite{MR3427047}
a blow-up theory is established for the minimizers of the functional
in~\eqref{9hiocdvjkb2toi3rg2fkfdjhsgddksjaqqd}, joined with an appropriate monotonicity formula,
which can be seen as a nonlocal counterpart of the classical
one in~\cite{weiss}
and which plays a crucial
role in establishing the homogeneity of the blow-up limit.

In addition, a classification result for
minimizing cones in the plane is given: namely,
if~$n=2$ and~$u$ is a continuous minimizer,
homogeneous of degree~$ s - \frac\sigma2>0$, then its free boundary
is a line.
\medskip

The one-phase case is then considered in~\cite{POIPOIPOI},
where a series of regularity results are obtained under
the additional assumption that the minimizer has a sign.
These results are summarized as follows:

\begin{theorem}\label{1fasee}
Let $(u,E)$ be a minimizer pair in~$\Omega\subset\R^n$
for the functional in~\eqref{9hiocdvjkb2toi3rg2fkfdjhsgddksjaqqd}, with~$u\ge0$.
Suppose that~$0\in\partial E$
and
\begin{equation}\label{LAMB01}
\int_{\R^n}\frac{|u(x)|}{1+|x|^{n+2s}}\,dx<+\infty.\end{equation}
Then:
\begin{itemize}
\item (Regularity of the minimizers). If~$s>\frac\sigma2$,
then~$u$ is locally H\"older continuous in~$\Omega$,
with H\"older exponent~$s-\frac\sigma2$.
\item (Density estimates for the free boundary). For
small~$r>0$, 
$$ \min\big\{ {\mathcal{L}}^n (E\cap B_r),\;
{\mathcal{L}}^n (E^c\cap B_r)\big\}\ge cr^n,$$
for some~$c>0$.
\end{itemize}
\end{theorem}

We think that it is an interesting problem to analyze the
case in which~$s\le\frac\sigma2$, and to
investigate
the two-phase problem in the setting of Theorem~\ref{1fasee}.
As a matter of fact, when the minimizer changes sign,
the two phases interact and produce new energy contribution
(and this phenomenon is new with respect to the classical case).
Also, the one-phase case of Theorem~\ref{1fasee} provides
the additional advantage that the condition~$u\ge0$ in~$E$
is automatically warranted by the sign of~$u$,
hence, there is a ``conceptual
freedom'' of choosing ``what $E:=\Omega_+(u)$
means'' in presence of ``fat zero level sets of~$u$''
(see again~\cite{POIPOIPOI}
for details).
\medskip

We also point out that condition~\eqref{LAMB01} is quite standard
in the fractional Laplace framework, since it is the ``natural''
one which ensures that the integrodifferential operator
in~\eqref{DEF:0123rugrfhijfoo98utytg} is convergent at infinity.
Nevertheless, recently a new theory for ``divergent'' fractional
Laplacians have been introduced in~\cite{DICDIVDIV}
and it would be interesting to reconsider several fractional
results also in the light of this new perspective.

\section{Nonlinear local plus nonlocal (or plus perimeter) free boundary problems}

In~\cite{2015arXiv151203043D}, energy functionals of the form
\begin{eqnarray}
\label{UNSK:1}
&&\int_\Omega |\nabla u(x)|^2\,dx+\Phi\Big( {\rm Per}_\sigma\big( 
\{u>0\},\Omega\big)\Big)\\
{\mbox{and }}\label{UNSK:2}
&&\int_\Omega |\nabla u(x)|^2\,dx+\Phi\Big( {\rm Per} \big(
\{u>0\},\Omega\big)\Big)
\end{eqnarray}
have been considered. Both the energy functionals
can be seen as a nonlinear energy superposition: the one in~\eqref{UNSK:1}
takes into account the fractional perimeter functional
with~$\sigma\in(0,1)$, which, in the
limit as~$\sigma\nearrow1$ recovers the classical perimeter in~\eqref{UNSK:2}.
Hence, to make the notation compact, in this section we will write~${\rm Per}_\sigma$
with~$\sigma\in(0,1]$ to denote the fractional perimeter when~$\sigma\in(0,1)$
and the classical perimeter when~$\sigma=1$
(also, the functionals in~\eqref{UNSK:1}
and~\eqref{UNSK:2} have to be suitably understood to take into consideration
competing pairs and boundary conditions, see~\cite{2015arXiv151203043D}
for details).\medskip

In this case, the Bernoulli condition along the free boundary
becomes both nonlinear and nonlocal. Namely, if~$u$ is
a minimizer for~\eqref{UNSK:1} or~\eqref{UNSK:2}
with~$u\ge0$ in~$E$,
along the free boundary
it holds that
\begin{equation}\label{234T1029dkjqO}
|\nabla u^+(x)|^2 -|\nabla u^-(x)|^2=
H_\sigma(x,E)\;
\Phi'\Big( {\rm Per}_\sigma \big(
\{u>0\},\Omega\big)\Big).
\end{equation}
Interestingly, condition~\eqref{234T1029dkjqO} 
reduces to~\eqref{BERN-SALSA}
when~$\sigma=1$ and~$\Phi(r)=r$, and to~\eqref{HGA:0}
when~$\sigma\in(0,1)$ and~$\Phi(r)=r$.\medskip

Once again, minimizers may vary from one scale to another,
due to the nonlinear effect in the energy functional
and if one does not impose any structural conditions
on the nonlinearity~$\Phi$ then singular free boundary may be developed also
in dimension~$2$, see Theorem~1.1 in~\cite{2015arXiv151203043D}.\medskip

On the other hand, for smooth and strictly
increasing nonlinearities,
several regularity results can be obtained, such as for instance the ones in
Corollary~1.4, Theorem~1.5 and Theorem~1.6 in\cite{2015arXiv151203043D}:

\begin{theorem}
Let $u$ be a minimizer for~\eqref{UNSK:1} or~\eqref{UNSK:2}
in~$\Omega\subset\R^n$
with~$u\ge0$ in~$E$ and~$0\in\partial E$.

Then:
\begin{itemize}
\item (Regularity of the minimizers). $u$ is locally H\"older continuous,
with H\"older exponent~$1-\frac\sigma2$.
\item (Density estimates for the free boundary). For
small~$r>0$, 
$$ \min\big\{ {\mathcal{L}}^n (E\cap B_r),\;
{\mathcal{L}}^n (E^c\cap B_r)\big\}\ge cr^n,$$
for some~$c>0$.
\end{itemize}\end{theorem}

As discussed in relation to Theorem~\ref{SIG:THM}, these results
do leave open many fundamental questions, such as the optimal regularity
for the minimizers and their free boundaries, which may be expected
to go beyond
the scaling properties of the energy.
In view of the considerations of this note,
in the future we also plan to provide a unified setting
to deal with general classes of possibly nonlinear and possibly nonlocal
free boundary problems.
\medskip

For completeness, we also mention that related but conceptually very
different types of nonlocal
free boundary problems have been recently considered in the literature
in the spirit of thin obstacle problems, in which a fractional equation is coupled
with a pointwise constraint of the solution, leading to a nonlocal type of variational
inequalities: for this, see e.g.~\cites{MR2367025, MR2511747, 2016arXiv160105843C}.

\vfill

\end{document}